\documentclass[10pt,a4paper]{article}

\usepackage{titletoc,amsmath,mathrsfs,amsfonts,amssymb}
\usepackage{txfonts}
\usepackage[pdftex]{hyperref}

\def\scr{\mathscr}
\newtheorem{thm}{\bf Theorem}[section]
\newtheorem{lem}[thm]{\bf Lemma}

\def\minus{\setminus}

\newcommand{\rf}[2]{[\ref{#1}; #2]}

\begin{document}

\thispagestyle{empty}

\begin{center}
{\bf\Large The first Dirichlet eigenvalue of birth-death process on tree}
\vskip.15in {Ling-Di WANG, Yu-Hui ZHANG}
\end{center}
\begin{center} wanglingdi@mail.bnu.edu.cn;\quad zhangyh@bnu.edu.cn\\
(Beijing Normal University, Beijing 100875, China)
 \end{center}

\bigskip
\noindent{\bf Abstract} {This paper investigates the birth-death (``B-D'' for short) process on tree with continuous time, emphasizing on estimating the principal eigenvalue (equivalently, the convergence rate) of the process with Dirichlet boundary at the unique root $0$. Three kinds of variational formulas for the eigenvalue are presented.
As an application, we obtain a criterion for positivity of the first eigenvalue for B-D processes on tree with one branch after some layer.
}

\section{ Introduction and main results}\label{c1}
This paper deals with   the first Dirichlet eigenvalue for B-D process on tree with the unique root $0$ as absorbing boundary. One may refer to \cite{Miclo, Ma} and the reference therein for more related works.  Our work is inspired by analogies research for B-D processes in \cite{r1,CWZ}, in which the principal eigenvalues in dimension one with kinds of boundary conditions were studied.
Let $T$ be a tree of at least two vertexes with the edge set $E$ (i.e., a connected graph without circle), such that the degree $d_i$ for each $i\in T$ is finite.  Let $|i|$ denote the  layer of $i$, and $i\sim j$ if $(i,j)\in E$. We call $j\in T$ a son (correspondingly, the father) of vertex $i\in T$ if $i\sim j$ and $|j|=|i|+1$ (correspondingly, $|j|=|i|-1$).
Consider a continuous time B-D process with $Q$-matrix such that $q_{ij}>0$ if and only if $i\sim j$. Then the corresponding operator is
 $$\Omega f(i)=\sum_{j\in J(i)}q_{i j} (f_j-f_i)+q_{i i^*}(f_{i^*}-f_i),\qquad i\in T,$$
 where  $J(i)$ is the set of  sons of $i$  and $i^*$ is the farther of $i$.
It is easy to obtain the unique symmetric measure $\mu$ on $T$:
$$\mu_0=1,\qquad\mu_k=\prod_{j\in\scr{P}(k)}\frac{q_{j^*j}}{q_{jj^*}},\quad k\in T\setminus\{0\},$$
where $\mathscr{P}(i)$  is the set of all the vertexes (the root $0$ is excluded) in the unique simple path from $i\in T\setminus\{0\}$ to the root.
If $(\lambda, g)$ with $g\neq0$ is a solution to ``eigenequation'':
\begin{equation}\label{eigen}\Omega g(i)=-\lambda g(i),\qquad i\in T\minus\{0\},\end{equation} then $\lambda$ is called  an ``eigenvalue'', and $g$ is called an ``eigenfunction'' of the eigenvalue $\lambda$. Note that the ``eigenvalue'' and ``eigenfunction'' used in this paper in a generalized sense rather than the standard ones since we do not require $g\in L^2(\mu)$. In this paper, we focus on estimating the principal Dirichlet eigenvalue $\lambda_0$ (i.e., the corresponding eigenfunction satisfies boundary condition $g_0=0$), which has the following classical variational formula:
\begin{equation}\label{f1}\lambda_0=\inf\{D(f): \mu(f^2)=1, f_0=0\},\end{equation}
where $\mu(f)=\sum_{k\in T\minus\{0\}}\mu_kf_k$ and
$$\aligned
D(f)=\sum_{i\in T\setminus\{0\}}\mu_iq_{ii^*}(f_i-f_{i^*})^2,\quad f\in\scr{D}(D)
\endaligned$$
with  $\scr{D}(D)=\{f: D(f)<\infty, f_0=0\}$.
Without loss of generality, we assume that the vertex $0$ has only one son throughout this paper (i.e., $|J(0)|=1$) and the  layer counting  begins from the son of the unique root $0$. Denote by  $N$ ($N\leqslant\infty$) the maximal layer of tree $T$ and $T_i$ ($i$ is included)  a subtree of tree $T$ with $i$  as root.  It is clear that $\lambda_0>0$ if $N<\infty$ (otherwise, $\Omega g(i)=0$. By letting $i\in E_N$ in \eqref{eigen}, we have $g_i=g_{i^*}$ for $i\in E_N$. By the induction, we have $g_i=g_0=0$ for $i\in T$, which is a contraction to $g\neq0$).
To state our results, we need some  notations as follows.  For $i\in T\setminus\{0\}$, define
$$\aligned
&I_i(f)=\frac{1}{\mu_iq_{ii^*}(f_i-f_{i^*})}\sum_{j\in T_i}\mu_jf_j \qquad \text{(single summation form)},\\
&I\!I_i(f)=\frac{1}{f_i}\sum_{k\in \scr{P}(i)}\frac{1}{\mu_kq_{kk^*}}\sum_{j\in T_k}\mu_jf_j \qquad \mbox{(double summation form)},\\
&R_i(w)=q_{ii^*}(1-w_i^{-1})+\sum_{j\in J(i)}q_{i j}(1-w_j)\qquad\mbox{(difference form)}.\endaligned$$
The forms of these operators defined above were initially introduced in \cite{Chen1996,Chen2001, r1} respectively for birth-death process in dimension one.  Shao and Mao in \cite{S-M} extended the operator with single summation form from line to tree, and obtained  the first operator defined above.
The domains of the three operators are defined respectively as follows:
$$\aligned
\mathscr{F}_I&=\{f: f_0=0, f_i>f_{i^*}\text{ for } i\in T\minus\{0\}\},\endaligned$$
$$\aligned\mathscr{F}_{I\!I}&=\big\{f: f_0=0, f>0 \text{ on } T\minus\{0\} \big\},\endaligned$$
$$\aligned\mathscr{W}&=\big\{w: w>1,  w_0=\infty\}.
\endaligned$$
These are used for the lower estimates of $\lambda_0$. For the upper bounds, some modifications are needed to avoid non-summable phenomenon, as shown below.
$$\aligned
\mathscr{{\widetilde F}}_{I}&\!=\!\big\{f\!>\!0:f_0\!=\!0, \exists 1\!\leqslant n\!<N+1 \text{ such that } f_i\!>\!f_{i^*} \text{ for } |i|\leqslant n, \text{ and }f_i\!=\!f_{i^*} \text{ for } |i|\!\geqslant\! n+1\big\},\endaligned$$
$$\aligned\mathscr{{\widetilde F}}_{I\!I}&=\{f>0: f_0=0, \exists1\leqslant n<N+1 \text{ such that } f_i=f_{i^*} \text{ for } |i|\geqslant n+1\},\endaligned$$
$$\aligned\mathscr{\widetilde{W}}&=\bigcup_{m:\,1\leqslant m<N+1}\bigg\{w: w_0\!=\!\infty, w_i>1\text{ and }\sum_{j\in J(i)}q_{ij}w_j<q_{i i^*}(1-w_i^{-1})+\sum_{j\in J(i)}q_{i j}\text{ for }
 \\
 &\hskip1.5cm |i|\leqslant m, \text{ and }w_i\!\!=\!1\text{ for } |i|\geqslant m+1\bigg\}.
\endaligned$$
Define $\widetilde R$ acting on $\widetilde{\scr{W}}$ as a modified form of $R$ by replacing $q_{i i^*}$ with $\mu_iq_{i i^*}\big/\sum_{j\in T_i}\mu_j$ in $R_i(w)$ when $|i|=m$,
where $m$ is the same one in $\widetilde{\scr{W}}$, when using approximating method, we also use $\widetilde R$ (at this time, $q_{ii^*}$ is replaced with $\tilde q_{ii^*}$ for each $i\in T$, see the arguments before Lemma \ref{L1}  and Step 4 in the proof of Theorem \ref{th1} below).
Here and in what follows, we adopt the usual convention $1/0=\infty$. The superscript ``$\widetilde{\quad}$'' means modified.

In  Theorem \ref{th1}  below, ``$\sup\,\inf$'' are
used for the lower bounds of $\lambda_0$, e.g., each test function $f\in\mathscr{F}_I$
produces a lower bound $\inf_{i\in T\minus\{0\}} I_i(f)^{-1}$, so this part
is called variational formula for lower estimate of $\lambda_0$. Dually, the
``$\inf\,\sup$'' are used for the upper estimates of $\lambda_0$. Among them, the ones expressed by  operator $R$ are easiest
to compute in practice, and the ones expressed by $I\!I$
are hardest to compute but provide better estimates.
Because of ``$\inf\,\sup$'', a localizing procedure is used
for the test function to avoid $I(f)\equiv \infty$ for instance,
which is removed out automatically for the ``$\sup\,\inf$'' part.
Define another set
$$\widetilde{\scr{F}}_{I\!I}'=\{f>0: fI\!I(f)\in L^2(\mu)\}.$$
Then we present our main results.
\begin{thm}\label{th1} The following variational formulas hold for $\lambda_{0}$ defined by \eqref{f1}.
\begin{itemize}
\item[$(1)$] Single summation  forms:
$$
\aligned \sup_{f\in\mathscr{F}_{I}}
\inf_{i\in T\setminus\{0\}}I_i(f)^{-1}=\lambda_{0}=
\inf_{f\in\mathscr{{\widetilde F}}_{I}}\sup_{i\in T\setminus\{0\}}I_i(f)^{-1},
\endaligned$$
\item[$(2)$] Double summation forms:
$$
\aligned
\sup_{f\in S(\mathscr{F})}\inf_{i\in T\setminus\{0\}}I\!I_i(f)^{-1}=\lambda_0&=\inf_{f\in S(\widetilde{\mathscr{F}})}\sup_{i\in T\setminus\{0\}}I\!I_i(f)^{-1}
\endaligned$$
with $S(\mathscr{F})=\mathscr{F}_{I\!I}$ or $\mathscr{F}_{I}$ and $S(\widetilde{\mathscr{F}})=\widetilde{\mathscr{F}}_{I\!I}$, or $\widetilde{\mathscr{F}}_{I}$,  or $\widetilde{\scr{F}}_{I\!I}'\cup\widetilde{\scr{F}}_{I\!I}$.
\item[$(3)$] Difference forms:
$$
\aligned
 \sup_{w\in\mathscr{W}}\inf_{i\in T\setminus\{0\}}R_i(w)=\lambda_{0}=\inf_{w\in\mathscr{\widetilde{W}}}\sup_{i\in T\setminus\{0\}}\widetilde R_i(w).
\endaligned$$
\end{itemize}
\end{thm}
We mention that the lower bounds of $\lambda_0$ in Theorem \ref{th1} (1) was known in \cite{S-M} as an inequality.  Liu et al. in \cite{L-M-W} extended the result in \cite{S-M}, obtained a  lower estimates  of $\lambda_0$ under some conditions.  In view of  the relation between the test functions of $R$, $I$ and $I\!I$ (they are all the analogies of eigenfunction, see arguments after Lemma \ref{L1} for details), it is not hard  to check that these estimates in Theorem ~\ref{th1}~can be sharp through the examples in \cite{S-M}, which illustrated that the lower estimates with single summation form was sharp.

Define $|A|=\text{number of elements in the set } A, $ $\mu(T_j):=\sum_{k\in T_j}\mu_k$, and
$$\varphi_j=\sum_{k\in \scr{P}(j)}\frac{1}{\mu_k q_{k k^*}}, \quad j\in T\minus\{0\}.$$
As applications of Theorem $\ref{th1}$ $(1)$ and $(2)$, we have the following theorem.
\begin{thm}\label{Basic}Let $\delta=\sup_{j\in T\minus\{0\}}\mu\big(T_j\big)\varphi_j$. Then
$$\delta^{-1}\geqslant\lambda_0\geqslant\bigg[\Big(2\sup_{i\in T\minus\{0\}}C_i\Big)\, \delta\bigg]^{-1},$$
where
$$C_i=1+|J(i)|+\sum_{s\in J(i)}\sum_{k\in T_s}\big(|J(k)|-1\big),\quad i\in T.$$
\end{thm}
The theorem effectively presents us the positive criterion of the first Dirichlet eigenvalue  for tree with one branch after some layer.
 For the degenerated case of the tree (only one branch), the results reduce to that of B-D process on half line studied in \cite{r1} (the ratio of the upper and lower bounds for the estimates of $\lambda_0$ is no more than $4$).  It is worthy to point out that the B-D process on tree with root as Dirichlet boundary can be a comparison with the B-D process on line with bilateral reflecting boundaries.  Let us have a look B-D process on line with reflecting boundaries. From \cite{r1}, we see that the eigenfunction of the first eigenvalue is strictly monotone with a unique zero.  If we treat the unique zero of the eigenfunction as root, then the B-D process on line is just B-D process on tree with two branches and the unique ``root'' as Dirichlet boundary (the intuition is pointed out by Professor Mao Y.H.). About the B-D process on line with reflecting boundaries, one may refer to \cite{Chen2013+B}.
\section{Proofs of the main results}
Define $E_m=\{i: |i|=m\}$, $T(n)=\cup_{m=0}^n E_m$ and
$$\tilde\lambda_0=\{D(f): \mu(f^2)=1, \exists1\leqslant n<N+1 \text{ such that } f_i=f_{i^*} \text{ for } |i|\geqslant n+1\}.$$
As will be seen  in Lemma \ref{L1} below, $\lambda_0=\tilde\lambda_0$ once $\sum_{k\in T}\mu_k<\infty$.
To this end,  define
$$\lambda_0^{(m)}=\inf\big\{D(f): \mu(f^2)=1, f_i=f_{i^*} \text{ for }\; |i|\geqslant m+1\big\},\qquad 1\leqslant m<N+1.$$
There is an explanation for $\lambda_0^{(m)}$ \big(see \rf{r1}{Section 4, Page 427}\big): let
$$\tilde\mu_i=\mu_i,\quad \tilde q_{i j}=q_{i j}\qquad\text{ for }\quad |i|\leqslant m-1\text{ and } |j| \leqslant m-1;$$
$$\tilde\mu_i=\sum_{j\in T_i}\mu_j,\quad \tilde q_{i^* i}=q_{i^* i},
\quad\tilde q_{i i^*}=\mu_{i} q_{i i^*}\bigg/\sum_{j\in T_i}\mu_j \qquad\text{ for }\quad |i|=m.$$
Noticing $\mu_i q_{i i^*}=\tilde\mu_i\tilde q_{i i^*}$,  for $f$ with $f_i=f_{i^*}$ for $|i|\geqslant m+1$, we have
$$D(f)=\sum_{i\in T(m)\minus\{0\}}\tilde\mu_i \tilde q_{i i^*}(f_i-f_{i^*})^2=:\widetilde D(f),\qquad \mu(f^2)=\sum_{i\in T(m)}\tilde\mu_i f_i^2=:\tilde\mu(f^2).$$
 So the $Q$-matrix $\widetilde Q=(\tilde q_{i j}: i,j\in T(m))$ is  symmetric  with respect to $\{\tilde\mu_i\}_{i\in T(m)}$ and $\lambda_0^{(m)}$ is the first Dirichlet eigenvalue of the local Dirichlet form $\big(\widetilde D, \scr{D}(\widetilde D)\big)$ with state space $T(m)$.

  For simplicity, we use ``iff'' to denote ``if and only if''  and $\uparrow$(resp. $\downarrow$) to denote increasing and decreasing
  throughout the paper.
\begin{lem}\label{L1} Assume that $\sum_{k\in T}\mu_k<\infty$.
We have $\lambda_0=\tilde\lambda_0$ and $\lambda_0^{(n)}\downarrow \lambda_0$ as $n\to N$.
 \end{lem}
 {\noindent{\it Proof}}\quad By definition of $\lambda_0$, for any $\varepsilon>0$, there exists $f$ such that
${D(f)}\big/{\mu(f^2)}\leqslant \lambda_0+\varepsilon$.
Construct $f^{(n)}$ such that $f_i^{(n)}=f_i$ for $|i|\leqslant n$ and $f_i^{(n)}=f_{i^*}$ for $|i|\geqslant n+1$. Since $\sum_{k\in T}\mu_k<\infty$, we have
$$\aligned
&D(f^{(n)})=\sum_{i\in T\minus\{0\}}\mu_iq_{ii^*}\big(f_i-f_{i^*}\big)^2=\sum_{i\in T(n)\minus\{0\}}\mu_iq_{ii^*}\big(f_i-f_{i^*}\big)^2\uparrow D(f)\\
&\mu({f^{(n)}}^2)=\sum_{i\in T(n)\minus\{0\}}\mu_if_i^2+\sum_{i\in E_{n+1}}\mu(T_i)f_{i^*}^2\to \mu(f^2).\endaligned$$
By definitions of $\lambda_0$, $\tilde\lambda_0$ and $\lambda_0^{(n)}$, the required assertion holds. $\qquad\Box$

 This lemma presents us an approximating procedure, makes it possible that sometimes  we only need to show some assertion or property holds for finite tree even if $N=\infty$ (see Step 6 (b) and Step 8 in proofs of Theorem \ref{th1} below). The following lemma was known in \cite{S-M}, gives us a important property of eigenfunction $g$. The property provides us the basis for the choices  of those test functions sets of operators $I$, $I\!I$ and $R$.
\begin{lem}$(\rf{S-M}{Proposition 2.4})$ For B-D process on tree $T$ (may have infinite vertexes). If $(\lambda_0, g)$ is a solution to \eqref{eigen} with boundary condition $g_0=0$ and $g\in L^2(\mu)$ holds, then $g_i>g_{i^*}$  for each $i\in T\minus\{0\}$.
\label{SM}\end{lem}

Obviously, for B-D process on finite tree $T$ (a tree with maximal layer $N<\infty$), the eigenfunction $g$ of the first Dirichlet eigenvalue satisfies $g_i>g_{i^*}$ for every $i\in T$.
 Before moving further, we introduce a general equation and discuss the origin of operators. Consider
$$\text{Poisson equation}: \qquad\Omega g(i)=-f_i,\qquad i\in T\minus\{0\}.$$
By multiplying $\mu_i$ on  both sides of the equation and making summation with respect to $i\in T_k\cap T(n)$ for some $k\in T\minus\{0\}$ with $|k|\leqslant n$, it is easy to check that
\begin{equation}\label{f3-I}\sum_{j\in E_{n+1}\cap T_k}\mu_j q_{j j^*}(g_{j^*}-g_j)+\mu_k q_{k k^*}(g_k-g_{k^*})=\sum_{j\in T_k\cap T(n)}\mu_j f_j,\qquad |k|\leqslant n.\end{equation}
 If $\lim_{n\to N}\sum_{j\in E_{n+1}\cap T_k}\mu_j q_{jj^*}(g_{j^*}-g_j)=0$ (which is obvious for $N<\infty$),
 then we obtain the form of the operator $I$ by letting $n\to N$ and $f=\lambda g$ in \eqref{f3-I}.  Moreover, if $g_0=0$ (which is clear for the eigenfunction of Dirichlet eigenvalue $\lambda_0$), then
$$g_i=\sum_{k\in\scr{P}(i)}\frac{1}{\mu_kq_{k k^*}}\sum_{j\in T_k}\mu_jf_j.$$
This explains where the operator $I\!I$ comes from. Similarly, from the eigenequation \eqref{eigen}, we obtain the operator $R$ by letting $w_i=g_i/g_{i^*}$. The eigenequation is a ``bridge'' among these operators.
 Based on \cite{r1, CWZ} and taking full advantage of these relations
  we present the proofs of the main results.

{\noindent\it Proof of Theorem $\ref{th1}$}\quad
 We introduce the following circle arguments for  lower bounds of $\lambda_0$.
$$\aligned\lambda_{0}&\geqslant\sup_{f\in\mathscr{F}_{I\!I}}\inf_{i\in T\setminus\{0\}}I\!I_i(f)^{-1}=\sup_{f\in\mathscr{F}_{I}}\inf_{i\in T\setminus\{0\}}I\!I_i(f)^{-1}=\sup_{f\in\mathscr{F}_{I}}
\inf_{i\in T\setminus\{0\}}I_i(f)^{-1}\geqslant\sup_{w\in\mathscr{W}}\inf_{i\in T\setminus\{0\}}R_i(w)\geqslant\lambda_0.
\endaligned$$

Step 1\quad Prove that $\lambda_{0}\geqslant\sup_{f\in\mathscr{F}_{I\!I}}\inf_{i\in T\setminus\{0\}}I\!I_i(f)^{-1}$

For  positive sequence $\{h_i\}_{i\in T\setminus\{0\}}$ and  $g$ with $g_0=0$, $\mu(g^2)=1$, we have
$$\aligned
1&=\sum_{k\in T\setminus\{0\}}\mu_kg_k^2=\sum_{k\in T\setminus\{0\}}\mu_k\bigg(\sum_{i\in\scr{P}(k)}\big(g_i-g_{i^*}\big)\bigg)^2\quad(\text{since } g_0=0)\\
&\leqslant\sum_{k\in T\setminus\{0\}}\mu_k\sum_{j\in\scr{P}(k)}\frac{\mu_jq_{jj^*}}{h_j}\bigg(g_j-g_{j^*}\bigg)^2
\sum_{i\in\scr{P}(k)}\frac{h_i}{\mu_iq_{ii^*}}\quad(\text{by Cauchy's ineq.})\\
&=\sum_{j\in T\setminus\{0\}}\mu_jq_{jj^*}\big(g_j-g_{j^*}\big)^2\frac{1}{h_j}\sum_{k\in T_j}\mu_k\sum_{i\in\scr{P}(k)}\frac{h_i}{\mu_iq_{ii^*}}\\
&\qquad(\text{by exchanging the order of sums,  and } j\in\scr{P}(k) \text{ iff }k\in T_j).\endaligned$$
For every $f$ with $fI\!I(f)<\infty$, let $h_i=\sum_{k\in T_i}\mu_kf_k$. By the proportional property, we get
$$\aligned\mu(g^2)&\leqslant D(g)\sup_{j\in T\setminus\{0\}}\bigg(\sum_{k\in T_j}\mu_k\sum_{i\in\scr{P}(k)}\frac{h_i}{\mu_iq_{ii^*}}\bigg)\bigg/{\sum_{k\in T_j}\mu_kf_k}\leqslant D(g) \sup_{j\in T\setminus\{0\}}I\!I_j(f)\endaligned$$
By $\eqref{f1}$, we have
$$\lambda_0\geqslant \inf_{j\in T\setminus\{0\}}I\!I_j(f)^{-1},\qquad f\in\mathscr{F}_{I\!I},$$
 and the required assertion follows by making the supremum with respect to $f\in\mathscr{F}_{I\!I}$.

Step 2\quad Prove that
$$\sup_{f\in\mathscr{F}_{I\!I}}\inf_{i\in T\setminus\{0\}}I\!I_i(f)^{-1}=\sup_{f\in\mathscr{F}_{I}}\inf_{i\in T\setminus\{0\}}I\!I_i(f)^{-1}=\sup_{f\in\mathscr{F}_{I}}
\inf_{i\in T\setminus\{0\}}I_i(f)^{-1}.$$

(a)\quad Prove that
$$\sup_{f\in\mathscr{F}_{I\!I}}\inf_{i\in T\setminus\{0\}}I\!I_i(f)^{-1}\geqslant\sup_{f\in\mathscr{F}_{I}}\inf_{i\in T\setminus\{0\}}I\!I_i(f)^{-1}\geqslant\sup_{f\in\mathscr{F}_{I}}
\inf_{i\in T\setminus\{0\}}I_i(f)^{-1}.$$
 The first inequality is clear since $\scr{F}_{I}\subset\scr{F}_{I\!I}$. Replacing $f$ in the denominator  of $I\!I_j(f)$ with $\sum_{k\in\scr{P}(j)}(f_k-f_{k^*})$. Using the proportional property, for $f\in\scr{F}_I$,  we have
$$\aligned\sup_{j\in T\setminus\{0\}}I\!I_j(f)&=\sup_{j\in T\setminus\{0\}}\bigg[\bigg(\sum_{k\in \scr{P}(j)}\frac{1}{\mu_kq_{kk^*}}\sum_{i\in T_k}\mu_if_i\bigg)\bigg/ {\sum_{k\in\scr{P}(j)}\big(f_k-f_{k^*}\big)}\bigg]\leqslant \sup_{k\in T\setminus\{0\}}I_k(f).\endaligned$$
So the required assertion holds.

(b)\quad To prove the equality, it suffices to show that
$$\sup_{f\in\mathscr{F}_{I}}
\inf_{i\in T\setminus\{0\}}I_i(f)^{-1}\geqslant\sup_{f\in\mathscr{F}_{I\!I}}\inf_{i\in T\setminus\{0\}}I\!I_i(f)^{-1}.$$
 For  $f\in\scr{F}_{I\!I}$, without loss of generality, assume that $I\!I(f)<\infty$. Let $g=fI\!I(f)$. Then $g\in\scr{F}_I$ and
 $$g_i-g_{i^*}=\frac{1}{\mu_iq_{ii^*}}\sum_{j\in T_i}\mu_jf_j\geqslant \sum_{j\in T_i}\mu_jg_j\inf_{k\in T\minus\{0\}}\frac{f_k}{g_k} ,\qquad i\in T\minus\{0\},$$
 and then the required assertion follows immediately since $f\in\scr{F}_{I\!I}$ is arbitrary.

There is another choice to show the equality. By Lemma \ref{SM}, we see that the eigenfunction $g$ satisfies that $g_i>g_{i^*}$ for $i\in T\minus\{0\}$ provided $N<\infty$. So $g\in \scr{F}_I$ and
 $\lambda_0=I_i(g)^{-1}$ for $i\in T\minus\{0\}$ (\rf{S-M}{Lemma 2.3}). By making the infimum with respect to $i\in T\minus\{0\}$ first and then the supremum with respect to $f\in\scr{F}_I$, we have
 $\lambda_0\geqslant\sup_{f\in\scr{F}_I}\inf_{i\in T\minus\{0\}}I_i(f)^{-1}$.
 There is a small gap in the proof since the eigenfunction $g$ may not belong to $L^2$  in the case of $N=\infty$. However, one may avoid this by a standard approximating procedure (according to the approximating idea used in Step 4 below). Combining this with Step 1 above, the required assertion follows immediately.

Step 3\quad Prove that
$\sup_{f\in\scr{F}_{I\!I}}\inf_{i\in T\minus\{0\}}I\!I_i(f)^{-1}\geqslant\sup_{w\in\scr{W}}\inf_{i\in T\minus\{0\}}R_i(w)$.

We first change the form of $R_i(w)$. For $w\in\scr{W}$, let $u$ with $u_0=0$ be a positive function on $T\minus\{0\}$  such that $w_i=u_i/u_{i^*}$ for $i\in T\minus\{0\}$, i.e.,
$$u_i=\prod_{j\in{\scr P}(i)}w_j\quad\text{ for }\quad i\in T\minus\{0\},\qquad u_0=0.$$
 Then $u_i>u_{i^*}$ for $i\in T\minus\{0\}$ and
\begin{equation}\label{f2-cir}R_i(w)=\frac{1}{u_i}\bigg[\sum_{j\in J(i)}q_{i j}(u_i-u_j)+{q_{i i^*}}\big(u_i-u_{i^*}\big)\bigg]=-\frac{\Omega u(i)}{u_i}.\end{equation}
Now we turn to the main text.
For any fixed $w\in\scr{W}$, without loss of generality, assume that $R(w)>0$. Let $u$ be a function mentioned above such that $w_i=u_i/u_{i^*}$ and $f=uR(w)>0$. Then $f\in\scr{F}_{I\!I}$ and
$\Omega u(i)=-f_i$.
Since $u_i>u_{i^*}$, by $\eqref{f3-I}$, we have
$$\sum_{j\in T_k\cap T(n)}\mu_jf_j\leqslant \mu_k q_{k k^*}(u_k-u_{k^*})<\infty,\qquad |k|\leqslant n.$$
 So $f\in L^1(\mu)$ and
$$u_k-u_{k^*}\geqslant \frac{1}{\mu_k q_{k k^*}}\sum_{j\in T_k}\mu_j f_j$$
by letting $n\to N$. Moreover,
$$u_i\geqslant \sum_{k\in \scr{P}(i)}\frac{1}{\mu_k q_{k k^*}}\sum_{j\in T_k}\mu_j f_j.$$
Hence,
$$\inf_{i\in T\minus\{0\}}R_i(w)=\inf_{i\in T\minus\{0\}}\frac{f_i}{u_i}\leqslant \inf_{i\in T\minus\{0\}} I\!I_i(f)^{-1},\qquad i\in T\minus\{0\},$$
and the required assertion follows  by making the supremum with respect to $f\in\scr{F}_I$ first and then with respect to $w\in\scr{W}$.

Step 4 \quad Prove that $\sup_{w\in\scr{W}}\inf_{i\in T\minus\{0\}}R_i(w)\geqslant\lambda_0$.

We first prove  that $\sup_{w\in\scr{W}}\inf_{i\in T\minus\{0\}}R_i(w)\geqslant0$. Let $f\in L^1(\mu)$ be a positive function on $T\minus\{0\}$ and $h=fI\!I(f)$ on $T\minus\{0\}$, $h_0=0$. Then $$h_i-h_{i^*}=\frac{1}{\mu_iq_{i i^*}}\sum_{k\in T_i}\mu_kf_k.$$
Put $\bar w_i=h_i/h_{i^*}$ for $i\in T\minus\{0\}$. By calculation, we have
$$\aligned-\Omega h(i)&=q_{i i^*}(h_i-h_{i^*})-\sum_{j\in J(i)}q_{j j^*}(h_j-h_{j^*})\quad(\text{by } j\in J(i) \text{ iff } i=j^*)\\
&=\frac{1}{\mu_i}\sum_{k\in T_i}\mu_kf_k-\sum_{j\in J(i)}\frac{q_{j^* j}}{\mu_j q_{j j^*}}\sum_{k\in T_j}\mu_k f_k\\
&=\frac{1}{\mu_i}\bigg[\sum_{k\in T_i}\mu_kf_k-\sum_{j\in J(i)}\sum_{k\in T_j}\mu_k f_k\bigg]\qquad(\text{by symmetric})\\
&=f_i. \endaligned$$
So
$$R_i(\bar w) = -\frac{\Omega h(i)}{h_i}=\frac{f_i}{h_i}>0,\qquad i\in T\minus\{0\},$$
and the required assertion then follows by making the infimum with respect to $i\in T\minus\{0\}$ first and then the supremum with respect to $w\in\scr{W}$.

By lemma \ref{SM}, if $\lambda_0>0$ and the maximal layer of the tree $N<\infty$, then  the eigenfunction $g$ satisfies $g_i>g_{i^*}$ for every $i\in T\minus\{0\}$.
Let $\bar w_i=g_i/g_{i^*}$. Then $\bar w\in\scr{W}$ and
$$R_i(\bar w)=\frac{-\Omega g(i)}{g_i}=\lambda_0,\qquad i\in T\minus\{0\}.$$
So the assertion holds for $N<\infty$. If $N=\infty$, then a approximating procedure is used.  Let $m\in \mathbb{N}^+$ and $1\leqslant m\leqslant N$. Then  $\lambda_0^{(m)}\downarrow\lambda_0$ as  $m\to N$ by Lemma \ref{L1}. Noticing the explanation of $\lambda_0^{(m)}$ at the beginning of this section and the assertion we have just  showed for $N<\infty$, we have
$$\lambda_0^{(m)}=\sup_{w\in\scr{W}(m)}\inf_{i\in T(m)\minus\{0\}}\widetilde R_i(w),$$
where $\scr{W}(m)=\{w: w_i>1, i\in T(m), w_0=\infty\}$, $\widetilde R$ is a modified form of $R$ by replacing $q_{ii^*}$ with $\tilde q_{ii^*}$.
By definition of supremum, for any fixed $\varepsilon>0$, there exists $\bar w\in\scr{W}$ such that
$$\lambda_0^{(m)}<\inf_{i\in T(m)\minus\{0\}}\widetilde R_i(\bar w)+\varepsilon\leqslant \inf_{i\in T(m-1)\minus\{0\}}\widetilde R_i(\bar w)+\varepsilon.$$
Extend $\bar w$ to $T$ by setting $\bar w_i=\bar w_{i^*}$ for $|i|>m$. Noticing $q_{ii^*}=\tilde q_{ii^*}$ for $|i|<m$, we have $\widetilde R_i(\bar w)=R_i(\bar w)$ for $|i|<m$.  Since
$\inf_{i\in T(m-1)\minus\{0\}}R_i(\bar w)\to\inf_{i\in T\minus\{0\}}R_i(\bar w)$ as $m\to \infty$, the required assertion follows by letting $m\to \infty$.

We adopt the following circle to prove the upper bounds of $\lambda_0$.
$$
\aligned
\lambda_0&\leqslant \inf_{f\in\widetilde{\scr{F}}_{I\!I}'\cup\widetilde{\scr{F}}_{I\!I}}\sup_{i\in T\minus\{0\}}I\!I_i(f)^{-1}\\
&\leqslant \inf_{f\in\widetilde{\scr{F}}_{I\!I}}\sup_{i\in T\minus\{0\}}I\!I_i(f)^{-1}=\inf_{f\in\widetilde{\scr{F}}_{I}
}\sup_{i\in T\minus\{0\}}I\!I_i(f)^{-1}=\inf_{f\in\widetilde{\scr{F}}_{I}
}\sup_{i\in T\minus\{0\}}I_i(f)^{-1}\\
&\leqslant \inf_{w\in\widetilde{\scr{W}}}\sup_{i\in T\minus\{0\}}
\widetilde{R_i}(w)\leqslant \lambda_0.
\endaligned$$
The second inequality above is clear, then we prove the remainders.

Step 5\quad Prove that $\lambda_0\leqslant \inf_{f\in\widetilde{\scr{F}}_{I\!I}'\cup\widetilde{\scr{F}}_{I\!I}}\sup_{i\in T\minus\{0\}}I\!I_i(f)^{-1}$.

For $f\in \widetilde{\scr{F}}_{I\!I}$, there exists $n\in E$ such that $f_i=f_{i^*}$ for $|i|\geqslant n+1$. Let $g_i=f_iI\!I_i(f)$ for $|i|\leqslant n$ and $g_i=g_{i^*}$ for $|i|\geqslant n+1$. Then $g\in L^2(\mu)$ and
$$g_i-g_{i^*}=\frac{1}{\mu_iq_{ii^*}}\sum_{j\in T_i}\mu_jf_j \mathbf{1}_{\{i: |i|\leqslant n\}}.$$
Inserting this term into $D(g)$, we have
$$\aligned D(g)&=\sum_{j\in T\minus\{0\}}(g_j-g_{j^*})
\sum_{k\in T_j}\mu_kf_k \mathbf{1}_{\{j: |j|\leqslant n\}}\\
&=\sum_{k\in T\minus\{0\}}\mu_kf_k\sum_{j\in\scr{P}(k)}\mathbf{1}_{\{j: |j|\leqslant n\}}(g_j-g_{j^*})\qquad(\text{since } k\in T_j \text{ iff } j\in\scr{P}(k))\\
&=\sum_{k\in T\minus\{0\}}\mu_kf_kg_k\qquad(\text{ since }g_i=g_{i^*} \text{ for } |i|\geqslant n+1).\endaligned$$
Since $g\in L^2(\mu)$, we further obtain
$$D(g)\leqslant \sum_{k\in T\minus\{0\}}\mu_k g_k^2 \sup_{k\in  T\minus\{0\}} \frac{f_k}{g_k}\leqslant \mu(g^2) \sup_{k\in  T\minus\{0\}}I\!I_k(f)^{-1}.$$
Hence,
$$\lambda_0\leqslant\frac{D(g)}{\mu(g^2)}\leqslant\sup_{k\in  T\minus\{0\}}I\!I_k(f)^{-1}.$$
The inequality above also holds for $f\in\widetilde{\scr{F}}_{I\!I}'$ since the key in the proof above is $g=fI\!I(f)\in L^2(\mu)$, which holds naturally for $f\in\widetilde{\scr{F}}_{I\!I}'$. So the assertion follows by making the infimum with respect to $f\in \widetilde{\scr{F}}_{I\!I}\cup\widetilde{\scr{F}}_{I\!I}'$ on  both sides of the inequality above.

Step 6\quad Prove that $$\inf_{f\in\widetilde{\scr{F}}_{I\!I}}\sup_{i\in T\minus\{0\}}I\!I_i(f)^{-1}=\inf_{f\in\widetilde{\scr{F}}_{I}
}\sup_{i\in T\minus\{0\}}I\!I_i(f)^{-1}=\inf_{f\in\widetilde{\scr{F}}_{I}
}\sup_{i\in T\minus\{0\}}I_i(f)^{-1}.$$

(a)\quad We first prove that
$$\inf_{f\in\widetilde{\scr{F}}_{I\!I}}\sup_{i\in T\minus\{0\}}I\!I_i(f)^{-1}\leqslant\inf_{f\in\widetilde{\scr{F}}_{I}
}\sup_{i\in T\minus\{0\}}I\!I_i(f)^{-1}\leqslant\inf_{f\in\widetilde{\scr{F}}_{I}
}\sup_{i\in T\minus\{0\}}I_i(f)^{-1}.$$
Since $\widetilde{\scr{F}}_{I}
\subset\widetilde{\scr{F}}_{I\!I}
$, the first inequality is clear. For $f\in\widetilde{\scr{F}}_{I}
$, there exists $1\leqslant n<N+1$ such that $f_i=f_{i^*}$ for $|i|\geqslant n+1$ and $f_i>f_{i^*}$ for $|i|\leqslant n$.
Since $f_i=\sum_{k\in\scr{P}(i)}(f_k-f_{k^*})$ for $|i|\leqslant n$, inserting this term to the denominator of $I\!I(f)$ and using the proportional property, we have
$$\aligned\inf_{i\in T\minus\{0\}}I\!I(f)&=\inf_{i\in T(n)\minus\{0\}}I\!I(f)\geqslant\inf_{i\in T\minus\{0\}}I_i(f). \endaligned$$
and the required assertion holds since $f\in\widetilde{\scr{F}}_{I}$ is arbitrary.

(b)\quad Prove the equality.

For $f\in\widetilde{\scr{F}}_{I\!I}$, $\exists 1\leqslant n<N+1$ such that $f_i=f_{i^*}$ for $|i|\geqslant n+1$ and $f>0$. Let $g_i=f_iI\!I_i(f)$ for $0<|i|\leqslant n$, $g_0=0$ and $g_i=g_{i^*}$ for $|i|\geqslant n+1$. Then $g\in\widetilde{\scr{F}}_{I}$ and
$$g_i-g_{i^*}=\frac{1}{\mu_iq_{i i^*}}\sum_{k\in T_i}\mu_j f_j,\qquad |i|\leqslant n.$$
Moreover,
$$\mu_iq_{i i^*}(g_i-g_{i^*})\leqslant \sum_{j\in T_i}\mu_j g_j\sup_{j\in T_i}\frac{f_j}{g_j}=\sum_{j\in T_i}\mu_j g_j\sup_{j\in T_i}I\!I_i(f)^{-1},\qquad i\in T\minus\{0\}.$$
Hence,
$$\sup_{k\in T\minus\{0\}}I_k(g)^{-1}\leqslant\sup_{k\in T\minus\{0\}}I\!I_k(f)^{-1}.$$
Then the assertion follows by making the infimum with respect to $g\in\widetilde{\scr{F}}_{I}$ first and then the infimum with respect to $f\in\widetilde{\scr{F}}_{I\!I}$.

Alternatively, there is another method to prove the equality. Combining with the arguments in Step 5 and Step 6 (a), it suffices to show that
$$\inf_{f\in\widetilde{\scr{F}}_I}\sup_{k\in T\minus\{0\}}I_k(f)^{-1}\leqslant\lambda_0.$$
To see this, assume that $g$ is an eigenfunction corresponding to $\lambda_0^{(m)}$. Then $g_i>g_{i^*}$ for $i\in T(m)$. Extend $g$ to the whole space by letting $g_i=g_{i^*}$ for $|i|\geqslant m+1$. Then $g\in \widetilde{\scr{F}}_I$ and
$$\lambda_0^{(m)}=\sup_{k\in T(m)\minus\{0\}}I_k(g)^{-1}=\sup_{k\in T\minus\{0\}}I_k(g)^{-1}\geqslant \inf_{f\in\widetilde{\scr{F}}_I}\sup_{k\in T\minus\{0\}}I_k(f)^{-1}.$$
Noticing Lemma \ref{L1}, the required assertion then holds by letting $m\to\infty$.

Step 7\quad Prove that $\inf_{f\in\widetilde{\scr{F}}_{I\!I}}\sup_{i\in T\minus\{0\}}I\!I_i(f)^{-1}\leqslant\inf_{w\in\widetilde{\scr{W}}}\sup_{i\in T\minus\{0\}}
\widetilde{R_i}(w)$.

First, we change the form of $\widetilde R$. For $w\in \widetilde{\scr{W}}$ with $w_i=1$ for $|i|\geqslant m+1$, let $g$ be a positive function on $T\minus\{0\}$ with $g_0=0$  such that $w_i=g_i/g_{i^*}$. Then $g_i>g_{i^*}$ for $|i|\leqslant m$ and $g_i=g_{i^*}$ for $|i|\geqslant m+1$. Since
$$\sum_{j\in J(i)}q_{ij}w_j<q_{i i^*}(1-w_i^{-1})+\sum_{j\in J(i)}q_{i j}\quad\text{ for }\quad |i|\leqslant m,$$
 we have
$\widetilde R_i(w)={-\widetilde\Omega g(i)}\big/{g_i}>0$ for $|i|\leqslant m$ and $\widetilde R_i(w)=0$ for $|i|\geqslant m+1$,
where $\widetilde\Omega$ is a change form of $\Omega$ with $q_{i j}$ by replacing $q_{i j}$ with $\tilde q_{i j}$ for $|i|\leqslant m+1$ and $|j|\leqslant m+1$.

Now, we come back to the main assertion.  For $w\in \widetilde{\scr{W}}$ with $w_i=1$ for $|i|\geqslant m+1$, let $g$ be the function mentioned above and
\begin{eqnarray*}
f_{i}=
\left\{
\begin{array}{lll}
\sum_{j\in J(i)}q_{i j}(g_i-g_j)+q_{i i^*}(g_i-g_{i^*}),&|i|\leqslant m-1, \\
\tilde q_{i i^*}(g_i-g_{i^*}),&|i|=m,\\
f_{i^*},&|i|\geqslant m+1.
\end{array}
\right.
\end{eqnarray*}
Then $f_i=-\widetilde\Omega g(i)>0$ for $|i|\leqslant m$. By \eqref{f3-I}, we have
$$\sum_{j\in E_{m}\cap T_k}\mu_j q_{j j^*}(g_{j^*}-g_j)+\mu_k q_{k k^*}(g_k-g_{k^*})=\sum_{j\in T_k\cap T(m-1)}\mu_j f_j,\qquad |k|\leqslant m-1.$$
 Since
$$\mu_iq_{i i^*}(g_i-g_{i^*})=\sum_{j\in T_i}\mu_jf_i=\sum_{j\in T_i}\mu_jf_j,\qquad|i|=m,$$
we have
$$\sum_{j\in E_{m}\cap T_k}\mu_j q_{j j^*}(g_{j^*}-g_j)=\sum_{j\in E_m\cap T_k}\sum_{i\in T_j}\mu_if_i=\sum_{j\in \big(T\minus T(m-1)\big)\cap T_k}\mu_jf_j,\quad |k|\leqslant m.$$
Hence, for $0<|k|\leqslant m$, we obtain
$$\mu_kq_{k k^*}(g_k-g_{k^*})=\sum_{j\in T_k\cap T(m-1)}\mu_j f_j+\sum_{j\in E_{m}\cap T_k}\mu_j q_{j j^*}(g_j-g_{j^*})=\sum_{j\in T_k}\mu_jf_j.$$
Moreover,
$$g_i=\sum_{k\in\scr{P}(i)}\frac{1}{\mu_kq_{k k^*}}\sum_{j\in T_k}\mu_jf_j,\qquad 0<|i|\leqslant m,$$
and $\widetilde R_i(w)=f_i/g_i=I\!I_i(f)^{-1}$ for $0<|i|\leqslant m$.
Since $\widetilde R_i(w)=0$ and $f_i=f_{i^*}$ for $|i|\geqslant m+1$, we obtain
$$\sup_{i\in T\minus\{0\}}\widetilde R_i(w)=\sup_{i\in T\minus\{0\}}I\!I_i(f)^{-1}\geqslant\inf_{f\in\widetilde{\scr{F}}_{I\!I}}\sup_{i\in T\minus\{0\}}I\!I_i(f)^{-1},\qquad w\in\widetilde{\scr{W}},$$
and the required assertion holds.

Step 8\quad Prove that $\inf_{w\in\widetilde{\scr{W}}}\sup_{i\in T\minus\{0\}}\widetilde R_i(w)\leqslant\lambda_0$.

 Let $g$ with $g_0=0$ be an eigenfunction of local first eigenvalue $\lambda_0^{(m)}$ and extend $g$ to  $T\minus\{0\}$ by setting $g_i=g_{i^*}$ for  $|i|\geqslant m+1$. Put $w_i=g_i/g_{i^*}$ for $i\in T\minus\{0\}$. Then $w\in\widetilde{\scr{W}}$. Since $m<\infty$, we have $\widetilde R_i(w)=\lambda_0^{(m)}>0$ for $i\in T(m)\minus\{0\}$, and $\widetilde R_i(w)=0$ for $T\minus T(m)$.
 Therefore,
 $$\aligned
 \lambda_0^{(m)}&=\max_{i\in T\minus\{0\}}\widetilde R_i(w)\\
 &\geqslant \inf_{w\in\widetilde{\scr{W}}: w_i=1\text{ for }|i|\geqslant m+1}\max_{i\in T(m)\minus\{0\}}\widetilde R_i(w)\\
 &\geqslant\inf_{w\in\widetilde{\scr{W}}: \exists n\geqslant1 \text{ such that } w_i=1\text{ for } | i |\geqslant n+1}\max_{i\in T\minus\{0\}}\widetilde R_i(w)\\
 &\geqslant \inf_{w\in\widetilde{\scr{W}}}\max_{i\in T \minus\{0\}}\widetilde R_i(w).
 \endaligned$$
 The assertion then follows by letting $m\to N$.\qquad$\Box$

Define $T_{i, j}=T_i\cup T_j$.  Then $T_{J(i)}=\big\{k: s\in J(i) \text{ and } k\in T_s\big\}$.
Similarly, we have $J(T_i)=\{k: s\in T_i\text{ and } k\in J(s)\}$.
It is obvious that $J(T_i)=T_{J(i)}$. Without loss of generality, we adopt convention that $\mu\big(T_k\big)=0$  if $T_k=\phi$. The proof of Theorem \ref{Basic}, which is an application of  Theorem \ref{th1}, is presented as follows.

{\noindent\it Proof of Theorem $\ref{Basic}$}\quad First, we prove that $\lambda_0^{-1}\leqslant\big(2\sup_{i\in T\minus\{0\}}C_i\big)\, \delta$. It is easy to see that
$$\aligned\sum_{j\in T_i}\mu_jf_j&=\sum_{j\in T_i}f_j\bigg[\mu(T_j)-\sum_{k\in J(j)}\mu(T_k)\bigg]\\
&=\sum_{j\in T_i}\mu(T_j)f_j-\sum_{j\in T_i}\sum_{k\in J(j)}\mu(T_k)f_j\\
&=\sum_{j\in T_i}\mu(T_j)f_j-\sum_{k\in T_{J(j)}}\mu(T_k)f_{k^*}\quad(\text{since } J({T_i})=T_{J(i)}\text{ and }k\in J(j) \text{ iff } j=k^*)\\
&=\mu(T_i)f_i+\sum_{k\in T_{J(i)}}\mu(T_k)\big(f_k-f_{k^*}\big)\quad\big(\text{since }T_i=\{i\}\cup T_{J(i)}\big).\endaligned$$
Put $f_j=\sqrt{\varphi_j}$ for $j\in T$. Then
$$\aligned\sum_{j\in T_i}\mu_j\sqrt{\varphi_j}&=\mu(T_i)\sqrt{\varphi_i}+\sum_{k\in T_{J(i)}}\mu(T_k)\big(\sqrt{\varphi_k}-\sqrt{\varphi_{k^*}}\big)\\
&\leqslant \delta\bigg[{\varphi_i}^{-1/2}+\sum_{k\in T_{J(i)}}\frac{1}{\varphi_k}\bigg(\sqrt{\varphi_k}-\sqrt{\varphi_{k^*}}\bigg)\bigg].
\endaligned$$
Since $\varphi_k\geqslant\varphi_{k^*}$, we obtain
$$\aligned
\sum_{k\in T_{J(i)}}\frac{1}{\varphi_k}\bigg(\sqrt{\varphi_k}-\sqrt{\varphi_{k^*}}\bigg)&\leqslant
\sum_{k\in T_{J(i)}}\big(\varphi_{k^*}^{-1/2}-\varphi_k^{-1/2}\big),
\endaligned$$
Noticing that $T_{J(i)}=J(T_i)$ and $k\in J(j)$ if and only if $k^*=j$, we have
$$\aligned
\sum_{k\in T_{J(i)}}\varphi_{k^*}^{-1/2}= \sum_{k\in J(T_i)}\varphi_{k^*}^{-1/2} &=\sum_{j\in T_{i}}\sum_{k\in J(j)}\varphi_j^{-1/2}=\sum_{j\in T_{i}}|J(j)|\varphi_j^{-1/2}.
\endaligned$$
Inserting this term to the inequality above, it is easy to see that
$$\aligned
\sum_{k\in T_{J(i)}}\frac{1}{\varphi_k}\big(\sqrt{\varphi_k}-\sqrt{\varphi_{k^*}}\big)
&\leqslant
\sum_{j\in T_{i}}|J(j)|\varphi_j^{-1/2}-\sum_{k\in T_{J(i)}}\varphi_k^{-1/2}\\
&=|J(i)|\varphi_i^{-1/2}+\sum_{k\in T_{J(i)}}\bigg(|J(k)|-1\bigg)\varphi_k^{-1/2}\\
&\leqslant\bigg[|J(i)|+\sum_{k\in T_{J(i)}}\bigg(|J(k)|-1\bigg)\bigg]\varphi_i^{-1/2}\quad\big(\text{since }\varphi_k\geqslant\varphi_{k^*}\big).
\endaligned$$
Hence,
$$\aligned\sum_{j\in T_i}\mu_j\sqrt{\varphi_j}&\leqslant \bigg[1+|J(i)|+\sum_{s\in J(i)}\sum_{k\in T_s}\bigg(|J(k)|-1\bigg)\bigg]\delta\varphi_i^{-1/2}
=C_i\delta\varphi_i^{-1/2}.
\endaligned$$
Since
$$\frac{1}{\sqrt{\varphi_i}-\sqrt{\varphi_{i^*}}}
=\frac{1}{\varphi_i-\varphi_{i^*}}\bigg(\sqrt{\varphi_i}+\sqrt{\varphi_{i^*}}\bigg)=\mu_i q_{i i^*}\big(\sqrt{\varphi_i}+\sqrt{\varphi_i}\big),$$
we obtain
$$\aligned
I_i(\sqrt\varphi)&=\frac{1}{\mu_iq_{i i^*}\big(\sqrt{\varphi_i}-\sqrt{\varphi_{i^*}}\big)}\sum_{j\in T_i}\mu_j\sqrt{\varphi_j}\\
&\leqslant C_i\delta\varphi_i^{-1/2}\bigg(\sqrt{\varphi_i}+\sqrt{\varphi_{i^*}}\bigg)\\
&\leqslant 2C_i\delta\qquad(\text{since } \varphi_i\geqslant \varphi_{i^*}).
\endaligned$$
It is clear that $\sqrt{\varphi}\in\scr{F}_I$, by Theorem \ref{th1}\,(1), we have
$$\lambda_0^{-1}\leqslant\inf_{f\in\scr{F}_I}\sup_{i\in T\minus\{0\}}I_i(f)\leqslant\sup_{i\in T\minus\{0\}}I_i(\sqrt{\varphi})\leqslant\Big(2\sup_{i\in T\minus\{0\}}C_i\Big)\delta.$$

Now, we prove that $\lambda_0\leqslant\delta^{-1}$.
For $i_0\in T\minus\{0\}$, let $f$ be a function such that
\begin{equation*} f_i=
\begin{cases}
    \varphi_i &  \text{if }\; i\in \scr{P}(i_0), \\
    \varphi_{i_0} &  \text{if}\; i\in T_{i_0},\\
    0   &\;  \text{Others}.
 \end{cases}
 \end{equation*}
Then
$$\aligned\sum_{j\in T_i}\mu_j f_j&=\sum_{j\in T_i\cap T_{i_0}}\mu_j\varphi_{i_0}+\sum_{k\in T_i\cap\big(\scr{P}(i_0)\minus\{0\}\big)}\sum_{j\in T_k}\mu_j\varphi_k. \endaligned$$
Since $f_i-f_{i^*}=1\big/\big(\mu_iq_{i i^*}\big)$ for $i\in\scr{P}(i_0)$ and $f_i-f_{i^*}=0$ for $i\in T\minus\scr{P}(i_0)$.
we have
$$
\aligned
\lambda_0^{-1}&=\sup_{g\in\widetilde{\scr F}_I}\inf_{i\in T\minus\{0\}}I_i(g)\geqslant\inf_{i\in T\minus\{0\}}I_i(f)\\
&=\inf_{ i\in\scr{P}(i_0)}\bigg(\sum_{j\in T_i\cap T_{i_0}}\mu_j\varphi_{i_0}+\sum_{k\in T_i\cap\big(\scr{P}(i_0)\minus\{0\}\big)}\sum_{j\in T_k}\mu_j\varphi_k\bigg)\\
&=\mu(T_{i_0})\varphi_{i_0},\qquad i_0\in T\minus\{0\}.
\endaligned$$
The assertion follows by making supremum with respect to $i_0\in T\minus\{0\}$ on the both sides of the inequality above.$\qquad\Box$

\medskip
{\small
\noindent{\bf Acknowledgements}  The work is supported in part by NSFC (Grant No.11131003), SRFDP (Grant No. 20100003110005), the ``985'' project from the Ministry of Education in China and the Fundamental Research Funds for the Central Universities. The first author thanks
 Y.T. Ma for her kind help for improving the upper estimate of $\lambda_0$ in  Theorem \ref{Basic}.  The both authors also thank professor M.F. Chen, Y.H. Mao for their useful suggestions during this work.
}

\end{document}